\documentclass{amsart}
\usepackage{graphicx}
\usepackage{latexsym}
\usepackage{hyperref}
\usepackage{amsfonts}
\usepackage[all]{xy}
\usepackage{amssymb, mathrsfs, amsfonts, amsmath}
\usepackage{amsbsy}
\usepackage{amsfonts}
\newtheorem*{acknowledgement}{Acknowledgement}
\newtheorem{corollary}{Corollary}

\newtheorem{lemma}{Lemma}

\newtheorem{problem}{Problem}
\newtheorem{theorem}{Theorem}
\newtheorem{example}{Example}
\numberwithin{equation}{section}

\begin{document}
\title[almost Ricci soliton]{Compact almost Ricci solitons with constant\\ scalar curvature are gradient}
\author{A. Barros$^{1}$,\,\,R. Batista$^{2}$\,\,and\,\,E. Ribeiro Jr.$^{3}$}
\address{$^{1}$ Departamento de Matem\'{a}tica-Universidade Federal do Cear\'{a}\\
60455-760-Fortaleza-CE-BR} \email{abdenago@pq.cnpq.br}
\thanks{$^{1,2}$ Partially supported by CNPq-Brazil}

\address{$^{2}$Current: Departamento de Matem\'{a}tica-Universidade Federal do Cear\'{a},
60455-760-Fortaleza-CE-BR\\
Permanent: Departamento de Matem\'{a}tica-Universidade Federal do Piau\'i, 64049-550-Teresina-PI-BR} \email{linmarcolino@gmail.com}

\address{$^{3}$ Departamento de Matem\'{a}tica-Universidade Federal do Cear\'{a},
60455-760-Fortaleza-CE-BR} \email{ernani@mat.ufc.br}
\thanks{$^{3}$ Partially supported by Postdoctoral program of IMPA-Brazil and FUNCAP-Brazil}
\urladdr{http://www.mat.ufc.br}
\date{February 11, 2013}

\begin{abstract}
The aim of this note is to prove that any compact non-trivial almost Ricci soliton $\big(M^n,\,g,\,X,\,\lambda\big)$ with constant scalar curvature is isometric to a Euclidean sphere $\Bbb{S}^{n}$. As a consequence we obtain that every compact non-trivial almost Ricci soliton with constant scalar curvature is gradient. Moreover, the vector field $X$ decomposes as the sum of a Killing vector field $Y$ and the gradient of a suitable function.
\end{abstract}

\maketitle

\section{Introduction}The study of an almost Ricci soliton was introduced in a recent paper due to Pigola et al. \cite{prrs}, where essentially they modified the definition of Ricci solitons by adding the condition on  the parameter $\lambda$ to be a variable function. More precisely, we say that a Riemannian manifold $(M^n,\,g)$ is an almost Ricci soliton, if there exist a complete vector field $X$ and a smooth soliton function $\lambda:M^n \to \Bbb{R}$ satisfying
\begin{equation}
\label{eqfund1}
R_{ij}+\frac{1}{2}(X_{ij}+X_{ji})=\lambda g_{ij},
\end{equation}
where $R_{ij}$ and $X_{ij}+X_{ji}$ stand, respectively, for the Ricci tensor and the Lie derivative  $\mathcal{L}_{X}g$ in coordinates. We shall refer to this equation as the fundamental equation of an almost Ricci soliton $\big(M^n,\,g,\, X,\,\lambda\big)$. It will be called \emph{expanding}, \emph{steady} or \emph{shrinking}, respectively, if  $\lambda<0$, $\lambda=0$ or $\lambda>0$. Otherwise, it will be called \emph{indefinite}.
When the vector field $X$ is a gradient of a smooth function $f:M^n \to \Bbb{R}$ the manifold  will be called a gradient almost Ricci  soliton. In this case the preceding equation turns out
\begin{equation}
\label{eqfund2}
{R_{ij}+\nabla_{ij}^{2}f=\lambda g_{ij}},
\end{equation}
where $\nabla_{ij}^{2}f$ stands for the Hessian of $f$.

Moreover, when $X$ is a Killing vector field the almost Ricci soliton will be called trivial. Otherwise it will be a nontrivial almost Ricci soliton. We notice that when $n\ge3$ and $X$ is a Killing vector field an almost Ricci soliton will be a Ricci soliton, since in this case we have an Einstein manifold, from which we can apply Schur's lemma to deduce that $\lambda$ is constant.
The rigidity result contained in Theorem $1.3$ of \cite{prrs} indicates that almost Ricci solitons should reveal a reasonably broad generalization of the fruitful concept of classical soliton. In fact, we refer the reader to \cite{prrs} to see some of these changes.

In the direction to understand the geometry of almost Ricci soliton, Barros and Ribeiro Jr. proved in \cite{br2} that a compact gradient almost Ricci soliton with non-trivial conformal vector field is isometric to a Euclidean sphere. In the same paper they proved an integral formula for compact case, which was used to prove several rigidity results, for more details see \cite{br2}. In \cite{catino},  Catino proved that a locally conformally flat gradient almost Ricci soliton, around any regular point of $f$, is locally a warped product with $(n-1)$-dimensional fibers of constant sectional curvature.

\begin{example}\label{exp1} In the compact case a simple example appeared in \cite{br2}. It was built over the standard sphere $\big(\Bbb{S}^n,\,g_{0}\big)$ endowed with the conformal vector field $X=a^{\top}$, where $a$ is a fixed vector in $\Bbb{R}^{n+1}$ and   $a^{\top}$ stands for its orthogonal projection over $T\Bbb{S}^n$. We notice that $a^{\top}$ is the gradient of the height function $h_a$; for more details see the quoted paper and Theorem 2.3 in \cite{prrs}.
\end{example}
It is well-known that all compact $2$-dimensional Ricci solitons are trivial, see \cite{hamilton1}. However, the previous example gives that there exists a non-trivial compact $2$-dimensional almost Ricci soliton.
On the other hand, given a vector field $X$ on a compact oriented Riemannian manifold $M^n$ the Hodge-de Rham decomposition theorem, see e.g. \cite{Warner},  gives that we may decompose $X$ as a sum of a gradient of a function $h$ and a divergence-free vector field $Y$, i.e.
\begin{equation}
\label{hdr}
X=\nabla h+Y,
\end{equation}where $div Y=0,$ see \cite{br2}. For simplicity let us call $h$ the Hodge-de Rham potential.

Now we present a generalization to integral formulae obtained in Theorem $4$ of \cite{br2} for the gradient case.

\begin{theorem}
\label{thm1}Let $\big(M^n,g,X,\lambda\big)$ be a compact  oriented almost Ricci soliton. If $S$ and $dV_{g}$ stand for the scalar curvature and the Riemannian volume form of $M^{n},$ respectively, then we have:
\begin{enumerate}
\item $\int_{M}|Ric-\frac{S}{n}g|^{2}dV_{g} =\frac{n-2}{2n}\int_{M}\langle\nabla S,X\rangle dV_{g}.$
\item $\int_{M}|Ric-\frac{S}{n}g|^{2} dV_{g}=\frac{n-2}{2n}\int_{M}\langle\nabla S,\nabla h\rangle dV_{g}.$
\end{enumerate}
\end{theorem}

As a consequence of the previous theorem we deduce a strong characterization to any compact almost Ricci soliton with constant scalar curvature, where $X$ is not necessarily the gradient of a potential function $f$. More precisely, we have the following corollary.

\begin{corollary}
\label{cor1thm1}
Let $\big(M^n,g,X,\lambda\big),\, n\ge3,$ be a  non-trivial compact oriented almost Ricci soliton. Then, the following facts are equivalent:
\begin{enumerate}
\item $M^{n}$ is isometric to a Euclidean sphere $\Bbb{S}^{n}.$
\item $\int_{M}\langle\nabla S,X\rangle dV_{g}=0.$
\item $\mathcal{L}_{X}S=0$, where $\mathcal{L}$ denotes Lie derivative.
\item $S$ is constant.
\item $M^n$ is homogeneous.
\end{enumerate}
\end{corollary}

We highlight that if a compact manifold has constant scalar curvature $S \leq 0,$ then every conformal vector field is Killing; see e.g. Theorem 6 in \cite{obata}. Therefore, the fourth assumption is relevant when the scalar curvature is positive. Actually,  Corollary  \ref{cor1thm1}  is used to give an answer to the following problem proposed in \cite{prrs}:

\begin{problem}
\label{prob1}
Under which conditions a compact almost Ricci soliton is necessarily gradient?
\end{problem}
As a consequence of Corollary \ref{cor1thm1} we have the following result which answers the previous problem for dimensions bigger than two.
\begin{corollary}
\label{cor2thm1}
Every compact almost Ricci soliton with constant scalar curvature is gradient.
\end{corollary}

Now we invoke Hodge-de Rham decomposition (\ref{hdr}) to write
\begin{equation}
\label{eqn1}
\frac{1}{2}\mathcal{L}_{X}g=\nabla^2h+\frac{1}{2}\mathcal{L}_{Y}g.
\end{equation}
In order to answer positively to Problem \ref{prob1} we have to find conditions which imply that $Y$ is a Killing vector field. An affirmative answer to this question is given below.

\begin{corollary}
\label{cor3}
Let $\big(M^n,g,X,\lambda\big),n\ge3,$ be a non-trivial compact oriented almost Ricci soliton. Then, in the Hodge-de Rham decomposition, $Y$ is a Killing vector field on $M^n$ provided at least one of the following conditions is satisfied:
\begin{enumerate}
\item $X$ is a conformal vector field.
\item $M^n$ has constant scalar curvature.
\end{enumerate}
\end{corollary}

\section{Preliminaries}

In this section we shall present some preliminaries which will be used during the paper. First we remember that for a Riemannian manifold $(M^{n},g)$ it is well-known the following lemma, see \cite{chow}.

\begin{lemma}\label{lem1}
Let $(M^{n},g)$ be a Riemannian manifold and $X$ a vector field in $M$. Then
\begin{equation}\label{eq1}
X_{ijk}-X_{ikj}=X_{t}R_{tijk}
\end{equation}
\begin{equation}\label{eq2}
    X_{ijkl}-X_{ikjl}=R_{tijk}X_{tl}+R_{tijk,l}X_{t}
\end{equation}
\begin{equation}\label{eq3}
    X_{ijkl}-X_{ijlk}=R_{tikl}X_{tj}+R_{tjkl}X_{it}
\end{equation}

\end{lemma}
Other important properties concern with of the Ricci tensor and of the scalar curvature are given below

\begin{equation}\label{eq4}
    R_{ij,k}=R_{ji,k}
\end{equation}
\begin{equation}\label{eq5}
     R_{ij,k}-R_{ik,j}=-R_{tijk,t}
\end{equation}
\begin{equation}\label{eq6}
    R_{ij,kl}-R_{ij,lk}=R_{tikl}R_{tj}+R_{tjkl}R_{it}
\end{equation}
\begin{equation}\label{eq7}
    \frac{1}{2}S_{k}=R_{ki,i}=R_{ik,i}
\end{equation}

Inspired on ideas developed in \cite{mrm} we may use the previous lemma to obtain the following result.
\begin{lemma}\label{lem2}
Let $(M^{n},g,X,\lambda)$ be an almost Ricci soliton. Then, the following formulae hold:
\begin{equation}\label{eq9}
    S+X_{ii}=n\lambda
\end{equation}
\begin{equation}\label{eq11}
    R_{lj}X_{l}=-X_{jii}-(n-2)\lambda_{j}
\end{equation}
\begin{equation}\label{eq12}
    R_{ij,k}-R_{ik,j}=-\frac{1}{2}R_{lijk}X_{l}+\frac{1}{2}(X_{kij}-X_{jik})+\lambda_{k}g_{ij}-\lambda_{j}g_{ik}
\end{equation}
\begin{proof}
In order to obtain (\ref{eq9}) it is enough to contract equation (\ref{eqfund1}). For equation (\ref{eq11}) computing the trace of equation (\ref{eq1}) in $i$ and $k$ we obtain
\begin{equation*}
    X_{iji}-X_{iij}=R_{liji}X_{l}=R_{lj}X_{l}.
\end{equation*}
Next, using fundamental equation (\ref{eqfund1}) we have
\begin{eqnarray*}
  R_{ij,i} &=& -\frac{1}{2}(X_{iji}+X_{jii})+\lambda_{i}g_{ij} \\
   &=& -\frac{1}{2}(X_{iji}-X_{iij}+X_{iij}+X_{jii})+\lambda_{i}g_{ij} \\
   &=& -\frac{1}{2}R_{lj}X_{l}-\frac{1}{2}(X_{iij}+X_{jii})+\lambda_{i}g_{ij}.
\end{eqnarray*}
Hence, using the twice contracted second Bianchi identity (\ref{eq7}) we deduce
\begin{equation*}
    \frac{1}{2}S_{j}=-\frac{1}{2}R_{lj}X_{l}-\frac{1}{2}(X_{iij}+X_{jii})+\lambda_{i}g_{ij},
\end{equation*}which enables us obtain equation (\ref{eq11}) after comparing the previous expression with covariant derivative of (\ref{eq9}).

Now, we derive equation (\ref{eq12}). Indeed, taking covariant derivative of (\ref{eqfund1}) we deduce
\begin{equation*}
    R_{ij,k}+\frac{1}{2}(X_{ijk}+X_{jik})=\lambda_{k}g_{ij}
\end{equation*}
and
\begin{equation*}
    R_{ik,j}+\frac{1}{2}(X_{ikj}+X_{kij})=\lambda_{j}g_{ik}.
\end{equation*}
Now we compare these previous expressions and we use equation (\ref{eq1}) to obtain
\begin{eqnarray*}
  R_{ij,k}-R_{ik,j} &=& -\frac{1}{2}(X_{ijk}+X_{jik}-X_{ikj}-X_{kij})+\lambda_{k}g_{ij}-\lambda_{j}g_{ik} \\
   &=& -\frac{1}{2}R_{lijk}X_{l}+\frac{1}{2}(X_{kij}-X_{jik})+\lambda_{k}g_{ij}-\lambda_{j}g_{ik}, \\
\end{eqnarray*}which concludes the proof of the lemma.
\end{proof}
\end{lemma}

Before to announce the next result we recall that the $X$-Laplacian of some tensor $T_{ik}$ on a Riemannian manifold $M^n$ is given by $$  \Delta_{X} T_{ik}= \Delta T_{ik} - T_{ik,s}X_{s}$$ for any vector field $X\in\mathfrak{X}(M).$ The next lemma is the main result of this section, it will be used to prove Theorem \ref{thm1}.

\begin{lemma}\label{lem3}
For an almost Ricci soliton $(M^{n}, g, X,\lambda)$ we have
\begin{eqnarray}\label{prin}
  \Delta_{X} R_{ik} &=&2\lambda R_{ik}-2R_{ijks}R_{js}+\frac{1}{2}R_{is}(X_{sk}-X_{ks})+\frac{1}{2}R_{sk}(X_{si}-X_{is})  \\
   &+&(n-1)\lambda_{ik}+\lambda_{jj}g_{ki}-\lambda_{ij}g_{kj}.\nonumber
\end{eqnarray}

\begin{proof}
Using equation (\ref{eq12}) we obtain
\begin{equation*}
    R_{ki,j}-R_{kj,i}=\frac{1}{2}R_{lkji}X_{l}+\frac{1}{2}(X_{jki}-X_{ikj})+\lambda_{j}g_{ki}-\lambda_{i}g_{kj},
\end{equation*}
taking the covariant derivative of the previous identity we have
\begin{equation}\label{eq14}
  R_{ki,jt}-R_{kj,it} = \frac{1}{2}(R_{ijkl,t}X_{l}+R_{ijkl}X_{lt})
   +\frac{1}{2}(X_{jkit}-X_{ikjt})+\lambda_{jt}g_{ki}-\lambda_{it}g_{kj}.
\end{equation}
On the other hand, from (\ref{eq6}) we deduce
\begin{eqnarray}\label{eq}
  R_{jk,ij} &=& R_{jk,ji}+R_{sjij}R_{sk}+R_{skij}R_{js}\nonumber \\
   &=& R_{jk,ji}+R_{si}R_{sk}+R_{skij}R_{js}.
\end{eqnarray}
Next, since $\Delta R_{ik}=R_{ik,jj},$ comparing the previous expression with (\ref{eq14}) we obtain
\begin{equation}
\label{eqdeltalem3}
    \Delta R_{ik}=R_{jk,ij}+\frac{1}{2}(R_{ijkl,j}X_{l}+R_{ijkl}X_{lj})+\frac{1}{2}(X_{jkij}-X_{ikjj})+\lambda_{jj}g_{ki}-\lambda_{ij}g_{kj}.
\end{equation}
Moreover, by second Bianchi identity we have
\begin{eqnarray*}
  R_{ijkl,j}X_{l} &=& -R_{ijlj,k}X_{l}-R_{ijjk,l}X_{l} \\
   &=& -R_{il,k}X_{l}+R_{ik,l}X_{l}.
\end{eqnarray*}
Comparing with (\ref{eqdeltalem3}) and using equation (\ref{eq}) we have
\begin{eqnarray*}
  \Delta R_{ik} &=& R_{jk,ij}+\frac{1}{2}(R_{ik,l}-R_{il,k})X_{l}+\frac{1}{2}R_{ijkl}X_{lj} \\
   &+& \frac{1}{2}(X_{jkij}-X_{ikjj})+\lambda_{jj}g_{ki}-\lambda_{ij}g_{kj} \\
   &=& R_{jk,ji}+R_{si}R_{sk}+R_{skij}R_{js}+\frac{1}{2}(R_{ik,l}-R_{il,k})X_{l} \\
   &+&\frac{1}{2}R_{ijkl}X_{lj} +\frac{1}{2}(X_{jkij}-X_{ikjj})+\lambda_{jj}g_{ki}-\lambda_{ij}g_{kj},
\end{eqnarray*}
thus, using the twice contracted Bianchi identity given by (\ref{eq7}) and (\ref{eqfund1}) we obtain
\begin{eqnarray}
\label{eq2prooflem3}
  \Delta R_{ik} &=& \frac{1}{2}S_{ki}+R_{si}R_{sk}+R_{skij}R_{js}-\frac{1}{2}R_{skij}X_{sj}\nonumber \\
   &+&\frac{1}{2}(R_{ik,l}-R_{il,k})X_{l}+\frac{1}{2}(X_{jkij}-X_{ikjj})+\lambda_{jj}g_{ki}-\lambda_{ij}g_{kj}\nonumber \\
   &=&\frac{1}{2}S_{ki}+R_{si}R_{sk}+R_{skij}R_{js}-R_{skij}(-R_{sj}+\lambda g_{sj}-\frac{1}{2}X_{js})\nonumber \\
   &+&\frac{1}{2}(R_{ik,l}-R_{il,k})X_{l}+\frac{1}{2}(X_{jkij}-X_{ikjj})+\lambda_{jj}g_{ki}-\lambda_{ij}g_{kj}\nonumber \\
   &=&\frac{1}{2}S_{ki}+R_{si}R_{sk}+2R_{skij}R_{js}+\lambda R_{ik}+\frac{1}{2}R_{skij}X_{js}\nonumber \\
   &+&\frac{1}{2}(R_{ik,s}-R_{is,k})X_{s}+\frac{1}{2}(X_{jkij}-X_{ikjj})+\lambda_{jj}g_{ki}-\lambda_{ij}g_{kj}.
\end{eqnarray}
Next, we computate the following sum
\begin{equation}\label{eq15}
    Y=\frac{1}{2}S_{ki}+R_{sk}R_{si}-\frac{1}{2}R_{is,k}X_{s}+\frac{1}{2}X_{skis}.
\end{equation}
First, taking the twice covariant derivative in (\ref{eq9}) we obtain
\begin{equation*}
    \frac{1}{2}S_{ki}=-\frac{1}{2}X_{ssik}+\frac{n}{2}\lambda_{ik},
\end{equation*}
comparing with (\ref{eq15}), (\ref{eq2}), (\ref{eq3}) and (\ref{eq5}) we have
\begin{eqnarray*}
  Y &=&-\frac{1}{2}X_{sski}+\frac{n}{2}\lambda_{ik}+R_{sk}R_{si}-\frac{1}{2}R_{is,k}X_{s}+\frac{1}{2}X_{skis} \\
   &=& \frac{1}{2}(X_{skis}-X_{sski})+R_{sk}R_{si}-\frac{1}{2}R_{is,k}X_{s}+\frac{n}{2}\lambda_{ik} \\
   &=& \frac{1}{2}(X_{skis}-X_{sksi}+X_{sksi}-X_{sski})+R_{sk}R_{si}-\frac{1}{2}R_{is,k}X_{s}+\frac{n}{2}\lambda_{ik} \\
   &=& \frac{1}{2}(R_{tsis}X_{tk}+R_{tkis}X_{st}+R_{tsks}X_{ti}+R_{tsks,i}X_{t}) \\
   &+&R_{sk}R_{si}-\frac{1}{2}R_{is,k}X_{s}+\frac{n}{2}\lambda_{ik} \\
   &=& \frac{1}{2}(R_{ti}X_{tk}+R_{tkis}X_{st}+R_{tk}X_{ti})+\frac{1}{2}(R_{sk,i}-R_{si,k})X_{s} \\
   &+&R_{sk}R_{si}+\frac{n}{2}\lambda_{ik}  \\
   &=& \frac{1}{2}(R_{si}X_{sk}+R_{tkis}X_{st}+R_{sk}X_{si})-\frac{1}{2}R_{tski,t}X_{s}+R_{sk}R_{si}+\frac{n}{2}\lambda_{ik} \\
   &=& R_{si}(R_{sk}+\frac{1}{2}X_{sk})+\frac{1}{2}(R_{tkis}X_{st}+R_{sk}X_{si})-\frac{1}{2}R_{tski,t}X_{s}+\frac{n}{2}\lambda_{ik} \\
   &=& R_{si}(-\frac{1}{2}X_{ks}+\lambda g_{sk})+\frac{1}{2}(R_{tkis}X_{st}+R_{sk}X_{si})-\frac{1}{2}R_{tski,t}X_{s}+\frac{n}{2}\lambda_{ik}\\
   &=&-\frac{1}{2}R_{si}X_{ks}+\lambda R_{ik}+\frac{1}{2}R_{tkis}X_{st}+\frac{1}{2}R_{sk}X_{si}-\frac{1}{2}R_{tski,t}X_{s}+\frac{n}{2}\lambda_{ik}.
\end{eqnarray*}
Substituting the previous expression in (\ref{eq2prooflem3}) we deduce
\begin{eqnarray*}
   \Delta R_{ik}&=&-\frac{1}{2}R_{si}X_{ks}+2\lambda R_{ik}+\frac{1}{2}R_{tkis}X_{st}+\frac{1}{2}R_{sk}X_{si}\\
   &-&\frac{1}{2}R_{tski,t}X_{s}+2R_{skij}R_{sj}+\frac{1}{2}R_{skij}X_{js}+\frac{1}{2}R_{ik,s}X_{s}  \\
   &-&\frac{1}{2}X_{ikss}+\lambda_{jj}g_{ki}-\lambda_{ij}g_{kj}+\frac{n}{2}\lambda_{ik}.
\end{eqnarray*}
Now from (\ref{eq2}) and (\ref{eq3}) we have
\begin{eqnarray*}
  X_{ikss}-X_{issk} &=& X_{ikss}-X_{isks}+X_{isks}-X_{issk} \\
   &=& R_{tiks}X_{ts}+R_{tiks,s}X_{t}+R_{tiks}X_{ts}+R_{tsks}X_{it}.
\end{eqnarray*}
On the other hand, taking the covariant derivative in (\ref{eq11}) we have
\begin{equation*}
    X_{issk}=-R_{ti,k}X_{t}-R_{ti}X_{tk}-(n-2)\lambda_{ik}.
\end{equation*}
Thus,
\begin{equation*}
    X_{ikss}=-R_{ti,k}X_{t}-R_{ti}X_{tk}+R_{tiks}X_{ts}+R_{tiks,s}X_{t}+R_{tiks}X_{ts}+R_{tk}X_{it}-(n-2)\lambda_{ik}.
\end{equation*}
Finally, we may use the first Bianchi identity and (\ref{eq5}) to infer
\begin{eqnarray*}
  \Delta R_{ik} &=& 2\lambda R_{ik}-2R_{ijks}R_{js}-\frac{1}{2}R_{si}X_{ks}+\frac{1}{2}R_{tkis}X_{st} \\
   &+&\frac{1}{2}R_{sk}X_{si}-\frac{1}{2}R_{tski,t}X_{s}+\frac{1}{2}R_{skij}X_{js}+\frac{1}{2}R_{ik,s}X_{s} \\
   &-&\frac{1}{2}(-R_{it,k}X_{t}-R_{it}X_{tk}+R_{tiks}X_{ts}+R_{tiks,s}X_{t}+R_{tiks}X_{ts}+R_{tk}X_{it}) \\
   &+&\frac{(n-2)}{2}\lambda_{ik}+\lambda_{jj}g_{ki}-\lambda_{ij}g_{kj}+\frac{n}{2}\lambda_{ik}\\
   &=& 2\lambda R_{ik}-2R_{ijks}R_{js}-\frac{1}{2}R_{si}X_{ks}+\frac{1}{2}R_{tkis}X_{st}+\frac{1}{2}R_{sk}X_{si} \\
   &-&\frac{1}{2}R_{tski,t}X_{s}+\frac{1}{2}R_{skij}X_{js}+\frac{1}{2}R_{ik,s}X_{s}+\frac{1}{2}R_{is,k}X_{s}+\frac{1}{2}R_{is}X_{sk}  \\
   &-&\frac{1}{2}R_{tiks}X_{ts}-\frac{1}{2}R_{tiks,s}X_{t}-\frac{1}{2}R_{tiks}X_{ts}-\frac{1}{2}R_{sk}X_{is}  \\
   &+&(n-1)\lambda_{ik}+\lambda_{jj}g_{ki}-\lambda_{ij}g_{kj}\\
   &=&2\lambda R_{ik}-2R_{ijks}R_{js}+\frac{1}{2}R_{is}(X_{sk}-X_{ks})+\frac{1}{2}R_{sk}(X_{si}-X_{is})\\
   &+&\frac{1}{2}R_{ik,s}X_{s}+\frac{1}{2}R_{ik,s}X_{s}-\frac{1}{2}R_{tisk,t}X_{s}-\frac{1}{2}R_{tski,t}X_{s}\\
   &+&R_{tkis}X_{st}-R_{tiks}X_{ts}-\frac{1}{2}R_{tiks,s}X_{t}+(n-1)\lambda_{ik}+\lambda_{jj}g_{ki}-\lambda_{ij}g_{kj}\\
   &=&2\lambda R_{ik}-2R_{ijks}R_{js}+R_{ik,s}X_{s}+\frac{1}{2}R_{is}(X_{sk}-X_{ks})+\frac{1}{2}R_{sk}(X_{si}-X_{is})\\
   &+&R_{skit}X_{ts}-R_{tiks}X_{ts}-\frac{1}{2}R_{tisk,t}X_{s}-\frac{1}{2}R_{tski,t}X_{s}\\
   &-&\frac{1}{2}R_{tkis,t}X_{s}+(n-1)\lambda_{ik}+\lambda_{jj}g_{ki}-\lambda_{ij}g_{kj}\\
   &=&2\lambda R_{ik}-2R_{ijks}R_{js}+R_{ik,s}X_{s}+\frac{1}{2}R_{is}(X_{sk}-X_{ks})+\frac{1}{2}R_{sk}(X_{si}-X_{is})\\
   &+&(n-1)\lambda_{ik}+\lambda_{jj}g_{ki}-\lambda_{ij}g_{kj},
\end{eqnarray*}
which finishes the proof of the lemma.
\end{proof}
\end{lemma}
We point out that Lemma \ref{lem3}  extends to non gradient almost solitons similar formulas obtained in Lemma 3.3 of \cite{prrs} and Lemma 2.3 of \cite{mrm}.

\section{Proof of the results}

\subsection{Proof of Theorem \ref{thm1}}
\begin{proof}
First of all we compute the trace of identity (\ref{prin}) to obtain
\begin{equation}\label{eq16}
    \frac{1}{2}\Delta S-\frac{1}{2}\langle\nabla S,X\rangle=\lambda S-|Ric|^{2}+(n-1)\Delta\lambda.
\end{equation}
Now, using that $|Ric-\frac{S}{n}g|^{2}=|Ric|^{2}-\frac{S^{2}}{n},$ we infer
\begin{equation}
\label{eqndeltas}
\frac{1}{2}(\Delta S-\langle\nabla S,X\rangle)=-|Ric-\frac{S}{n}g|^{2}+\frac{S}{n}(n\lambda-S)+(n-1)\Delta\lambda.
\end{equation}
On integrating identity (\ref{eqndeltas}) and using  (\ref{eq9}) and the compactness of $M^{n}$ we arrive at
\begin{equation}
\label{eqn1thm1}
    \int_{M}|Ric-\frac{S}{n}g|^{2} dV_{g}=\frac{1}{2}\int_{M}\langle\nabla S,X\rangle dV_{g}+\frac{1}{n}\int_{M} S div X dV_{g}.
\end{equation}
Now we recall that for any vector field $Z$ on $M^n$ we have
\begin{equation}\label{eqn2thm1}
div(SZ)=Sdiv Z+\langle\nabla S,Z\rangle.
\end{equation} Whence, we deduce
\begin{equation}
\label{eqn3thm1}
    \int_{M}|Ric-\frac{S}{n}g|^{2} dV_{g}=\frac{n-2}{2n}\int_{M}\langle\nabla S,X\rangle dV_{g}=-\frac{n-2}{2n}\int_{M}SdivX dV_{g},
\end{equation}
which finishes the first statement. Proceeding, we notice that using once more identity (\ref{hdr}) we have from (\ref{eqn3thm1}) $\int_{M}|Ric-\frac{S}{n}g|^{2} dV_{g}=-\frac{n-2}{2n}\int_{M}S\Delta h dV_{g} $, which completes the proof of the theorem.
\end{proof}

\subsection{Proof of Corollary \ref{cor1thm1}}
\begin{proof}
We point out that any assumption of the corollary jointly with Theorem \ref{thm1} give $Ric=\frac{S}{n}g$. Whence we obtain $\frac{1}{2}\mathcal{L}_{X}g=(\lambda-\frac{S}{n})g,$ i.e., $X$ is a non-trivial conformal vector field. Now we may apply  a classical theorem due to Nagano and Yano in \cite{NY}  to conclude that $M^{n}$ is isometric to a Euclidean sphere (see also Theorem $2$ in \cite{br2}), which finishes the proof of the corollary.
\end{proof}

\subsection{Proof of Corollary \ref{cor2thm1}}

\begin{proof}
First we notice that from Corollary \ref{cor1thm1} we may assume that $M^n$ is isometric to the standard sphere of constant curvature $1$ and $X$ is a non-trivial conformal vector field.  In particular, using (\ref{eqfund1}) and the Hodge-de Rham decomposition $X=\nabla h +Y$ we obtain
\begin{equation}
\label{eqndeltah}\Delta h= n\lambda -n(n-1).
\end{equation} Now we invoke identities   (\ref{eqndeltas}) and (\ref{eqndeltah}) to obtain
\begin{equation}
\label{eqndeltahl}
\Delta (h+\lambda)=0.
\end{equation}
Whence we have $h=-\lambda +c$, where $c$ is a constant. Using once more (\ref{eqndeltah})
we deduce
\begin{equation}
\label{eqndeltah1}
\Delta h+nh=n(c-(n-1)).
\end{equation}
Therefore, up to additive constant, $h$ is a first eigenfunction of the Laplacian of $\Bbb{S}^{n}$. It is now standard to obtain that  $\nabla h$ is a conformal vector field and, according to (\ref{eqndeltah}), we have $$\frac{1}{2}\mathcal{L}_{\nabla h}g=\frac{\Delta h}{n}g=\big(\lambda-(n-1)\big)g.$$
Since $X$ is also conformal satisfying $$\frac{1}{2}\mathcal{L}_{X}g=\big(\lambda-(n-1)\big)g,$$ we may invoke Hodge-de Rham decomposition to deduce that $Y$ is a Killing vector field, as required. So, we complete the proof of the corollary.
\end{proof}

\subsection{Proof of Corollary \ref{cor3}}
\begin{proof} When $X$ is a conformal vector field it is well known that $\int_{M}\langle\nabla S,X\rangle dV_{g}=0,$ see e.g. \cite{bezin}. Hence, we deduce from Theorem \ref{thm1} that  $Ric=\frac{S}{n}g.$ Thus we have that its scalar curvature is constant and the result follows from the last part of the proof of Corollary \ref{cor2thm1}. The second assertion follows from the last argument and we complete the proof of the corollary.

\end{proof}

\begin{acknowledgement}The third author would like to thank the IMPA-Brazil for its support, where part of the work was started. He would like to extend his special thank to Professor F.C. Marques for very helpful conversations.  Finally, the authors want to thank the referees for their careful reading and helpful suggestions. \nonumber
\end{acknowledgement}

\end{document}